# Low power in-situ AI Calibration of a 3 Axial Magnetic Sensor

Roger Alimi[1], Elad Fisher[1], Amir Ivry[2], Alon Shavit[3], Eyal Weiss[1]

[1]Technology Division, Soreq NRC, Yavne 81800, Israel
[2]Andrew and Erna Viterbi Faculty of Electrical Engineering, Technion – Israel Institute of Technology, Haifa 3200003, Israel
[3]School of Computer Science and Engineering, Hebrew University of Jerusalem, Israel

**Magnetic surveys are conventionally performed by scanning a domain with a portable scalar magnetic sensor. Unfortunately, scalar magnetometers are expensive, power consuming and bulky. In many applications, calibrated vector magnetometers can be used to perform magnetic surveys. In recent years algorithms based on artificial intelligence (*AI*) achieve state-of-the-art results in many modern applications. In this work we investigate an *AI* algorithm for the classical scalar calibration of magnetometers. A simple, low cost method for performing a magnetic survey is presented. The method utilizes a low power consumption sensor with an *AI* calibration procedure that improves the common calibration methods and suggests an alternative to the conventional technology and algorithms. The setup of the survey system is optimized for quick deployment in-situ right before performing the magnetic survey. We present a calibration method based on a procedure of rotating the sensor in the natural earth magnetic field for an optimal time period. This technique can deal with a constant field offset and non-orthogonality issues and does not require any external reference. The calibration is done by finding an estimator that yields the calibration parameters and produces the best geometric fit to the sensor readings. A comprehensive model considering the physical, algorithmic and hardware properties of the magnetometer of the survey system is presented. The geometric ellipsoid fitting approach is parametrically tested. The calibration procedure reduced the root-mean-squared noise from the order of $10^4$ nT to less than 10 nT with variance lower than 1 nT in a complete 360º rotation in the natural earth magnetic field. In a realistic survey scheme the obtained calibration noise is suited to the environmental survey clutter. Implementing this scheme with a modern low power analog-to-digital convertor and micro-controller results in power consumption lower than 15 mW and calibration duration of few minutes.**

## I. INTRODUCTION

Magnetometers are used in surveys to map anomalies in the natural earth magnetic field [1]. Magnetic anomalies are caused by gradients in the magnetic permeability made by objects containing ferromagnetic mass [1], [2]. They may be used to map archeological sites [3], prehistoric sites [4], underground pipes [5], and even wrecks [3]. The location of the anomaly caused by ferromagnetic objects such as Unexploded Ordnance (UXOs) can be estimated by employing Magnetic Anomaly Detection (MAD) and localization schemes [3], [7].

Magnetic surveys are conventionally performed by scanning a domain with a portable scalar magnetic magnetometer. The scalar sensor measures the amplitude of the magnetic field and is not sensitive to the field bearing. Unfortunately, conventional scalar magnetometers are expensive, power consuming and bulky. Furthermore, most sensitive scalar magnetometers have "dead zones" orientation that renders them impractical in some applications.

A practical approach to overcome the shortcomings of scalar magnetometers is to use vector magnetometers. A 3-axis fluxgate magnetometer [1], [4], [5] is a low-cost, low-power and high resolution vector sensor [6] and is therefore a straight forward candidate for this application. Unfortunately, due to mechanical and sensitivity inaccuracies between the axes, the 3 sensors' axes are not precisely orthogonal [7], [8]. These imperfections lead to very significant magnitude errors when the sensor spins in 360 degrees. Typical values of non-orthogonality are less than 1º, offset error of less than 10 nT and scaling error of less than 0.5%. Therefore, it is impossible to use the fluxgate sensor for our purposes, without calibration.

Theoretically the offset of the sensor may be measurable inside a perfectly shielded environment. However, since the residual field inside a practical shield is directional and non-uniform, it interferes with the offset of each axis. In order to compensate for the residual field, a calibration must be performed within the shield itself.

Calibration procedures are commonly used. They can be classified into two main approaches: the scalar approach sensor in which the reference is a scalar precise magnetometer [7], and the vector approach in which the reference values come from a pre-corrected or ideal vector magnetometer [9].

Mathematical solutions were proposed by Marklund et al. for calibration of satellites magnetometers during magnetic mapping missions [10]. The proposed solutions assume that the intensity of the magnetic field where the magnetometer is rotated is precisely known. Other solutions, aimed for spacecraft equipped with a fluxgate sensor, use a scalar reference or a known constant field [14] to enable real time calibration during typical spacecraft mission mode operations [11]. Merayo et al report a peak to peak of 0.5 nT and a variance of about 0.1 nT, while Gavazzi et al. achieved a 0.3 nT standard deviation in different environment [12]. However, those methods use an expensive scalar magnetometer, while compelling the user to use an additional external device. Furthermore, the power consumption of the suggested solutions is relatively high.

Other methods suggest a calibration which is based on a procedure of moving the sensor inside a constant (unknown) magnetic field. A comprehensive model which generalizes the problem has been suggested, and verified that the determination of all calibration parameters from a motion inside a constant magnetic field is possible [13].

The geometric approach calibration is a scalar calibration technique that was presented in several articles. This technique is capable of calibrating with a constant field offset. The calibration is done by finding an estimator which yields the calibration parameters, which produce the best geometric fit to the sensor readings. This technique does not require any external reference [14]. Some calibration algorithms based on a geometric approach have been published [15], and some were even demonstrated recently in some applications like calibrating a compass mechanism of underwater gliders [16].

In this work, we present a low-power, low-noise and low-cost magnetometer optimized to allow both conventional vector measurements and calibrated scalar measurements. We experimentally compare between two approaches: a classical geometric one and a more advances computational algorithm based on an artificial intelligence (AI) scheme. Both of them allow in-situ calibration of the vector magnetometer.

Thus, a low cost, low power consumption and small size fluxgate three-axial magnetometer in magnetic surveys can be achieved. A comprehensive model considering the physical, algorithmic and hardware properties of the magnetometer of the survey system is presented. The geometric ellipsoid fitting approach [17]–[19] is parametrically tested. An experimental system was built and tested using the suggested in situ low power calibration algorithm.

## II. MAGNETIC CALIBRATION

### A. Calibration Goal

Our goal is to present a low power in situ calibration process and to compare two different calibration concepts: a standard deterministic calculation and one AI approach, namely a Neural Network (NN) scheme. We do this by optimizing a low power sensor, sampling rate and a calibration algorithm to enable a calibration process in-situ. We use two parameters to model an error in a measurement. The first is the difference between the highest to the lowest magnetic magnitude, maximum peak to peak error (PTP). The second parameter is the norm of the magnetic field variance. Both parameters are calculated from the measured data during a given measurement scenario.

The low-power and in-situ optimization are important for highly portable magnetic survey applications. In order to perform a magnetic survey, the measurement system signal-to-noise ratio best be greater than 1 [1]. The most dominant sources of noise in land surveys are clutter in the ground and inhomogeneity of the ground ferrous composition giving rise to PTP clutter of 1-4 nT. This clutter is a source of noise as the survey is not aimed to find metal ferromagnetic scrap. Our goal is to reach a PTP calibration error of the same magnitude as of the clutter noise where the portable magnetometer is carried manually across the domain as stable as may be expected. Additionally, in order to allow quick deployment, the calibration procedure should consume less than 5 minutes and the power consumption should be lower than 15 mW.

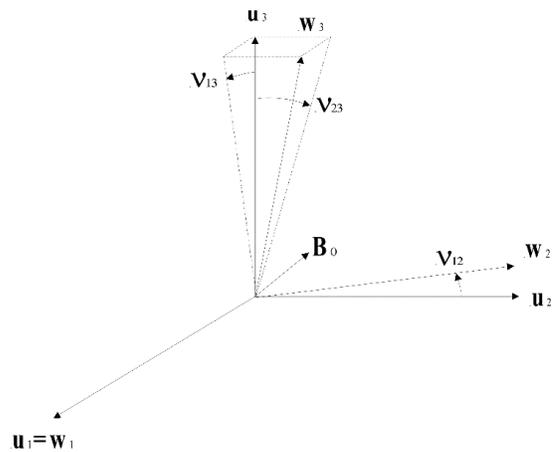

Fig. 1. The different coordinate systems of a fluxgate magnetometer. The orthogonal system is defined by the axes u1, u2 and u3, while the non-orthogonal system is defined by the axes W1, W2 and W3. The orthogonality errors are modeled by the angles, $v_{12}, v_{23}$ and $v_{13}$.

### B. Magnetometer

In this work, we focus on the simplest possible scanning configuration, in which the scanning process is done manually by walking over an area. The sensor is being held by hand during all the mapping process. Hence, while mapping an area, the sensor changes its orientation.

The simplest way to conduct a magnetic mapping is by using scalar magnetic sensors, such as Geometrics' G-823A cesium magnetometer. These types of sensors have a PTP of roughly 0.3 nT over 360° rotation and 0.2 nT during a stationary measurement, which renders them "over-spec" for magnetic mapping purposes. However, scalar sensors have "dead zones" in some sensor orientations, and an expensive tag price in the range of 10,000-40,000$ as well as a relatively high-power consumption of about 12W. Another alternative is using a fluxgate sensor. This is a low-priced (200-500$) vector sensor, which its output consists of the magnetic field components.

As illustrated in Fig. 1, the coordinate system defined by the axes *W1*, *W2* and *W3*, is not orthogonal. The non-orthogonality errors between the axes can be modeled by the angles, $v_{12}$, $v_{23}$ and $v_{13}$.

## III. CALIBRATION APPROACHES

In this work we present two different algorithms to perform the calibration and compare their performances. We apply the calibration procedures on the same measured data set. Moreover, we have also generated a synthetic, augmented dataset using a simulation setup, which allows a deeper analysis of the AI algorithm performance.

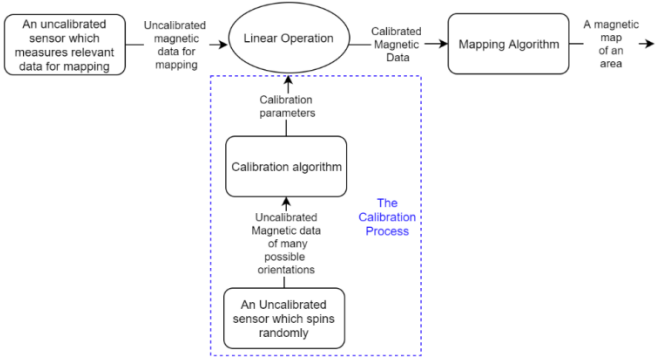

Fig. 2. Calibration data flow.

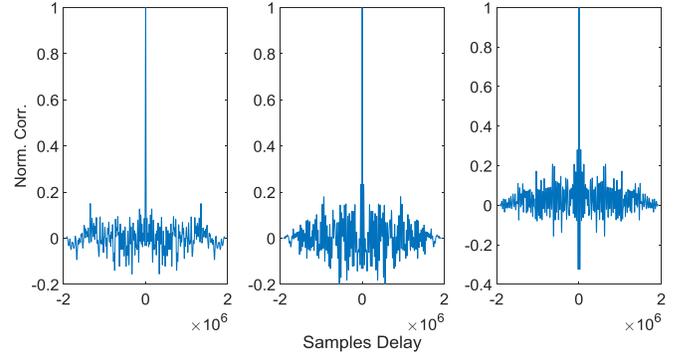

Fig. 3. Normalized cross-correlation between the measurements of each 3-axial component of $L_e$ and $L_s$. From left to right: $x$ component (RMSE = 100 nT), $y$ component (RMSE = 77 nT), $z$ component (RMSE = 89 nT).

### A. Geometrical approach

This solution is based on a linear transformation, which is applied on the sensor's raw data. This transformation solves the problem of the non-orthogonal axes by converting each vector components of the magnetic field at a given time $\left(B_x(t), B_y(t), B_z(t)\right)$ into a new corrected data. Fig. 2 presents a block diagram which demonstrates the data flow in our solution. A scalar magnetometer, which measures a constant magnetic field, would show a constant magnitude $|B|$, irrespective of the magnetometer's orientation. Thus, each measurement value has the same distance from the origin:

$$|B| = \sqrt{B_X(t)^2 + B_Y(t)^2 + B_Z(t)^2}, \; \forall t \in \mathbb{R}. \tag{1}$$

Hence, when plotting the magnetic data of different sensor orientations, a perfect sphere should be formed. However, a realistic fluxgate magnetometer has some non-orthogonality between its axes. Readings of a non-ideal magnetometer are always on an ellipsoid manifold [20]. Therefore, by using a linear transformation, it is possible to cast the calibration problem to a unified transformation parameterized by rotation ($R$), scaling ($S$) and offset ($b$). This approach produces a linear problem that may be solved by using conventional linear methods. The rotation and scaling are due to the non-orthogonality between the axes and temperature changes, while the offset (the ellipsoid's displacement from the origin) is caused due to constant ferromagnetic variations formed by the magnetometer core itself.

The solution consists of an algorithm that performs a linear transformation of the magnetometer readings. The algorithm fits all the points into an ellipsoid manifold centered on the origin. After applying the transformation, all the data measurement points lie at the same distance from the origin. Since this distance is defined by $|B|$, the result is a constant magnetic magnitude, irrespective of the sensor's orientation [21], [22].

There are several approaches concerning how to transform an ellipsoid into a sphere while minimizing the fitting error. One method is the "Ellipsoid fit" presented in [22]. This algorithm describes the ellipsoid using nine parameters: semi-axis $(a, b, c)$, Euler angles $(\varphi, \theta, \Psi)$ (repressing the successive axis rotations) and coordinates of the center $(X_0, Y_0, Z_0)$. The fitting is done by performing a linear convex optimization on the raw data $h_r$ to minimize the error of the fit, by minimizing the sum of squared algebraic error distances (*RMS*).

According to the algorithm output, the correction matrix ($M$) and the offset vector ($b$), which best transform the ellipsoid into an origin-centered sphere are calculated. $M$ is calculated by using $R'$ (de-rotation) and $S^{-1}$ (de-scaling). The corrected data ($h_c$) is then given by:

$$h_c = M(h_r - b), \; s.t. \, M = S^{-1}R'. \tag{2}$$

A constant magnetic field must be used during all time of calibration. For this purpose, we take the median value of all the earth's magnetic field measurements (about 2 million data points) which can be considered as a quite stable and precise estimation of the "true" ideal value.

### B. Artificial Intelligence Approach

In this section, we present an additional calibration algorithm. In contrast to the algorithm presented in the previous section, this method will map each 3D experimental data point on a 3D ideal components that lie on a sphere the radius of which is the median value of the measured norms of the field. Let us notate the set of noisy un-calibrated measurements as $L_e$ and the set of their perfectly calibrated counterparts as $L_s$. While each 3-axial component in $L_e$ lays on an ellipsoid, each component in $L_s$ lays on a sphere that is centered at the origin with radius that is equal to the earth's magnetic field $B$. Therefore, we perform a pre-processing step to find the magnetic measurements that comprise this sphere, i.e., to generate the set $L_s$. An optimization problem is solved, in which every noisy measurement in $L_e$ is mapped to its closest counterpart on the sphere, in terms of the Euclidean distance. Formally, for a given noisy coordinate $(x_e, y_e, z_e) \in L_e$, the matched clean coordinate $(x_s, y_s, z_s) \in L_s$ is obtained

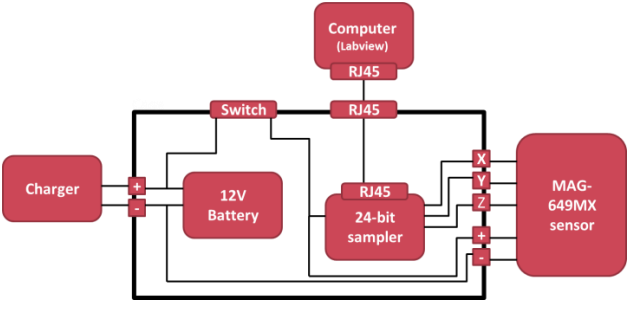

Fig. 4. Hardware configuration scheme.

by solving the following non-convex optimization problem:

$$(x_s, y_s, z_s) = \underset{(x,y,z): x^2+y^2+z^2=B^2}{\operatorname{argmin}} J(x,y,z,x_e,y_e,z_e) \quad (3)$$

where $J = (x - x_e)^2 + (y - y_e)^2 + (z - z_e)^2$ and $B$ stands for the earth's magnetic field. To ensure a non-trivial solution, the cross-correlation between $L_e$ and $L_s$ is given in Fig. 3, for each of the 3-axial components. To emphasize the existing difference between these curves, the root mean squared error (*RMSE*) is also generated.

In recent years, neural networks have re-gained popularity due to an increase in the amount of available data and computational capabilities [29], [30]. In this study, we attempt to perform calibration by learning the mapping from the set $L_e$ to the set $L_s$ with an artificial neural network. The network comprises two hidden layers of three nodes each. Thus, the output layer, which holds the perfectly calibrated 3-axial magnetic measurements in $L_s$, is given by the multiplication of the input layer by two $3 \times 3$ matrices. These matrices are notated as $S_{net}^{-1}$ and $R_{net}$. Let $(a_e, b_e, c_e)$ and $(a_s, b_s, c_s)$ represent the parameter sets of the original ellipsoid and the target sphere, respectively. Then, the network learns the following mapping:

$$\Psi: \frac{x_e^2}{a_e^2} + \frac{y_e^2}{b_e^2} + \frac{z_e^2}{c_e^2} \rightarrow \frac{(x_s - a_s)^2 + (y_s - b_s)^2 + (z_s - c_s)^2}{B^2}. \quad (4)$$

Namely, the mapping $\Psi$ has an intrinsic linear nature. By taking this fact and the noisy measurements that lay on the ellipsoid into account, the activation chosen for the network is linear and the objective function is the mean-squared error (MSE). By integration of the drop-out approach on 2 neurons every 30 epochs, the network is effectively trained to avoid over-fitting [23]. Thus, the mapping between the ellipsoid and the sphere is modelled with maximal generalization ability, and the statistical properties of the noise do not affect the model. Given an unseen coordinate $x_e^{test} \in L_e$, the following mapping is applied to generate $x_s^{test} \in L_s$:

$$x_s^{test} = \Psi(x_e^{test}): x_e^{test} S_{net}^{-1} R_{net} + b_{net}, \quad (5)$$

where $b_{net}$ stands for a bias component.

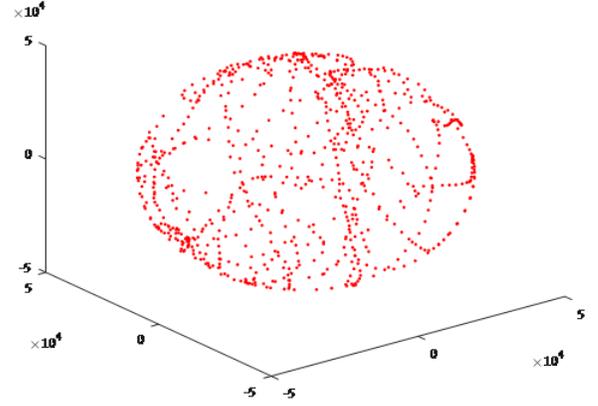

Fig. 5. Points which were created during the sensor's spinning for the calibration process. The axes are the x, y and z components of the magnetic field. As expected, the low-density ellipsoid shape is observed.

IV. EXPERIMENTAL

Through this section both authentic (in-situ) and synthetic datasets are considered in order to accomplish valid, real-world applicable results on one hand and demonstrate the potential of the suggested *AI* method on a variety of setups, on the other.

*A. Apparatus*

Bartingtons' Mag648 and Mag 649 low power fluxgates were used [24]. They were digitized by a 24-bit NI-9239 DAQ [25]. The sampled signal was then transferred to a computer by LabVIEW® data processing application. The NI sampling unit was used to simulate high performance sampling scheme in order to determine the optimal sampling requirements. A practical sampling unit such as TI 1292 [26] was tested as well for low power optimization. The hardware configuration we used is presented in [1].

*B. Measurements and simulation setups*

To collect genuine calibration data, the sensor is rotated by hand for approximately 2 minutes in the natural undisturbed earth field. This rotation is performed in a ferromagnetic free environment. The collected data is used as an input for the calibration algorithm. The algorithm's output is a correction matrix and an offset, which are specific for each magnetometer. After applying the matrix and the offset on the sensor's raw data, we get the calibrated data. This calibration process can be done in the field in real time "in-situ" and allow immediate use the magnetometer as a scalar survey tool.

To generate synthetic data, a MATLAB® simulation was built. This simulation generates uniformly sampled data points given a set of pre-determined parameters such as the sensors' internal noise level, data acquisition time, analog and digital frequencies, etc. Once again, this data serves as input for the calibration algorithms, which ultimately produce correction matrices to be applied to any new raw data.

*C. Parameters optimization*

During the calibration process the sensor is rotated at a rate of 1 revolution per second with rotation modes at higher frequencies. This revolution rate might be too high for the conventional low sampling frequency [1] and may result in "smearing" of the measurement and thus decreasing the accuracy of the calibration. Thus, the sampling rate and the averaging windows must be optimized in order to maintain signal integrity.

TABLE 1. THE EFFECT OF THE SAMPLE FREQUENCY AND BANDWIDTH ON PTP AND VARIANCE USING THE GEOMETRICAL APPROACH.

| Configuration | PTP (nT) | Var (nT²) |
|---|---|---|
| Mag649 (250 Hz) | 15 | 7 |
| Mag649 (3000Hz) | 15 | 3 |
| Mag648 (250 Hz) | 38 | 23 |
| Mag648 (1000 Hz) | 29 | 23 |

TABLE 2. CALIBRATION PERFORMANCE OF GEOMETRIC AND NEURAL NETWORK COMPARED TO UNCALIBRATED MEASUREMENT: MAGNITUDE MAX PTP AND VARIANCE. MOVING AVERAGE WINDOW WITH 0.1 SEC WIDTH IS APPLIED (MAG649 AT 3 KHZ).

| | PTP ($nT$) | Var ($nT^2$) |
|---|---|---|
| Not Calibrated | 24,055 | 652 |
| Geometric | 15 | 3 |
| Neural network | 8 | 4 |

We used a faster sampling hardware to optimize the sampling rate and averaging window, as shown in Table 1 shows the main findings of this study. A sampling frequency of 3kHz provides the best results.

Another important parameter is the bandwidth of the sensor. This is because the high rate of rotation during the calibration introduced artificial signals at frequencies beyond the bandwidth of the sensor. To test this, we have used 2 variants of the same sensor by Bartington instruments: Mag648 with bandwidth of 5 Hz and Mag649 with bandwidth of 1 kHz. Table 1 shows the main findings of this study. A sampling frequency of 3kHz provides the best results.

It appears that Mag649 with the larger bandwidth and frequency response band has lower PTP error and variance compared to Mag648. Furthermore, increasing the sampling rate decreases the calibration error (using the geometric approach algorithm for comparison) in both narrowband and wide band sensors. Table 1 shows the main findings of this study. A sampling frequency of 3kHz provides the best results.

Averaging of random noise decreases the noise. However, averaging over a long period of time averages not only random noise but also the signal which is changing during the measurement period. As a result, the signal integrity may decrease. We have tested different averaging windows to optimize the tradeoff between signal integrity and noise decrease. The best value was taken to be a moving average window of 0.1 sec.

*D. Experiments*

Several experiments are conducted to analyze the performance of both geometric and *AI* methods in-situ, as well as on synthetic simulated data.

Initially, the geometric and *AI* methods are employed on the authentic data measurements, and their performance is given in Table 2.

Next, a set of synthetic measurements is generated using a dedicated simulation. This allows the analysis of the *AI* approach without the limitations the dictates from the real-world setups. It should be highlighted that the **uniform** sampling done in the simulation differs from the **random** spinning of the sensor, done in-situ. The same coverage **amount** even if it is random in both cases, will give better results for a true uniformly distributed process. This can be achieved easily in simulation but it is more complicated to achieve for manual in-situ setup. A mechanical mechanism with uniform sampling data points pre-calculated is presented in [27].

To acquire synthetic measurements, two sensor noise levels $\sigma_n = 0.1$ and 0.3nT, and coverage percentages, ranging from 10% up to 100% is used. The definition of coverage percentage in this study refers to the area of the sphere covered by the measurements. Let the sphere be parameterized by the set $\{B, \theta, \phi\}$, where $B$ is the magnetic field and $0 \leq \theta \leq 2\pi, 0 \leq \phi \leq \pi$. It is known that an area differential value on the spheres' surface is proportional to $\frac{\Delta\phi\Delta\theta}{B}$. In this study we make use of this and define a coverage of 100% when $\frac{\Delta\phi\Delta\theta}{B} \leq 10^{-4}$ for all bald areas on the surface of the sphere. Here, bald area is an area in which no measurement lays. The measurements then undergo the *AI* calibration process, and their performance output are extracted.

The coverage percentage is highly correlated with the systems' end-to-end calibration process computational time. The more time for calibration the higher the accuracy, given the environment is noise free. This allows the examination of the effectiveness of uniform sampling data acquisition in comparison its real-world counterpart. In our case, we have physically rotated the sensor for 2 minutes. The end-to-end calibration process using the simulation was done on a Core-i7-7820HQ CPU 64-bit operating system, x64 based processor.

V. RESULTS

*A. Measurements*

The calibration algorithm transforms the measurements from an ellipsoid to a sphere. The measured ellipsoid can be seen in Fig. 5. An effective calibration algorithm transforms 3 non-orthogonal axial sensor output to scalar data. The quality of the calibration algorithm may be measured by the data peak to peak value compared to not-calibrated output, as well as the variance of the calibrated data with respect to $B$.

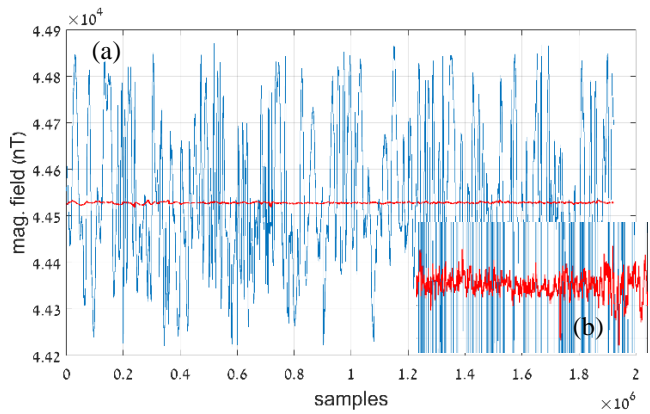

Fig. 6. Magnetic magnitude over 360 degrees' sensor spin, before (a) and after applying the geometrical calibration [central line zoomed in (b)]. Sensor used is Mag649 with sample rate of 3kHz.

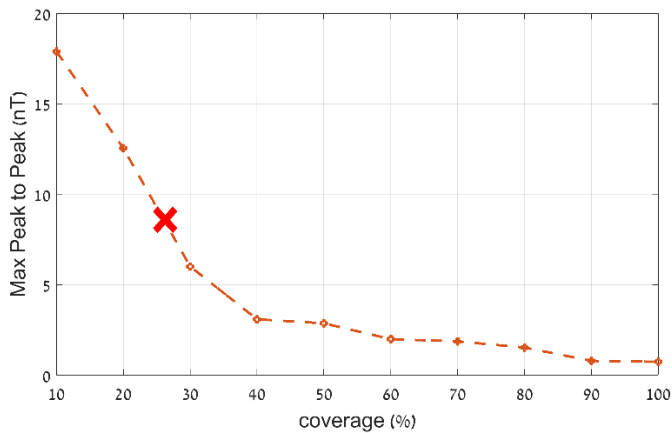

Fig. 7. Peak to peak calibration performance of AI approach for various sphere coverages and two sensor internal noise values. The X sign denotes the experimental peak to peak. It lies at 26% coverage of the full sphere surface.

### B. Conventional Geometrical approach

A graphical demonstration of the input data and the output of the calibration algorithm can be seen in Fig. 6 (raw data vs. calibration results). Here we present the results for the selected configuration of the system. The configuration yielding the lowest PTP noise and variance was using a wide band sensor (Mag649) with fast sampling rate of 3000 Hz and averaging over 0.1s window. Implementing this scheme with a low power analog-to-digital convertor (*ADC*) and micro-controller results in power consumption lower than 15 mW, and calibration duration of approximately 2 minutes. The best result is a maximum peak to peak value of 15 nT with a variance of 3 $nT^2$.

### C. Artificial Intelligence Approach

To demonstrate the performance of the artificial intelligence method, two sets of calculations are conducted, as extended in the previous section. In the in-situ experiment, the *AI* approach shows a better performance over its geometric oriented opponent: a maximum peak to peak of 8 nT has been obtained (instead of 15nT for the geometric) and a maximum variance of 4 nT (similar to the one obtained with the previous method).

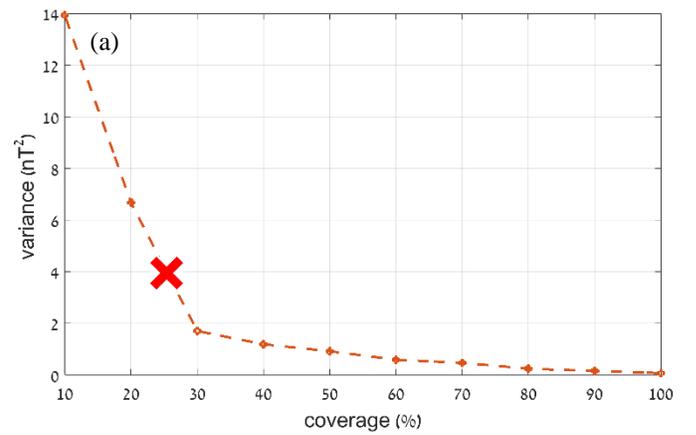

Fig. 8. Variance calibration performance of *AI* approach for various sphere coverages and two sensor internal noise values. The X notation denotes the experimental peak to peak. It lies at 26% coverage of the full sphere surface.

When considering the simulation related results, a larger set of comparisons can be made. For instance, in Figs. 7 and 8 it is clearly shown that the *AI* approach highly improves when given a dataset that covers the sphere even by 70%, **as long as it is done uniformly**. By extending it to 60% coverage, the mean error value approaches to the same order of magnitude that the natural environmental clutter noise [1].

## VI. Discussion

Several conclusions can be drawn based on these results. First, let us address the results achieved with the experimental in-situ measurements. As previously stated, the *NN* learns the linear mapping from the ellipsoid to the sphere. Therefore, the performance of this calibration model is limited in several ways when it comes to non-uniform sampling. First and foremost, a new sample can be mapped by the network to an interval of values with uniform probability, where these values may be distant from $B$. This is a direct outcome of the linear fit the network applies on its training data. Second, because the model effectively maps an ellipsoid onto a sphere, the noisy test measurements lead to a non-optimal variance. Also, it can be concluded that using sparse training set, which contains many uncovered regions in the ellipsoid, is limited in its performance and leads to mediocre results considering the vast number of in-situ recorded samples, ranging up to almost 2 million.

When converging to full coverage, one can achieve a performance that approaches the known barrier of the sensors' internal noise. That and more, the presented system can handle authentic data in real-time, under the constraint of uniform sampling. e.g., by allowing a delay of 60 minutes, the achieved error approaches the system' hardware limit with a negligible error variance of less than $0.5\ nT^2$ (see Fig. 8). However, in magnetic surveys, this is not important because the environmental clutter is significantly larger than the sensors' internal noise.

These outcomes highly emphasize the gap in effectiveness between uniform and random data sampling during acquisition. i.e., using the 2 million experimental samples lead to the same performance achieved by using merely 26% uniform sphere coverage (see Fig. 7 and Fig. 8). This essentially means that if the authentic data acquisition is done uniformly, it can save almost 80% of its computational resources and achieve similar performance.

From the measurements and theoretical projections, a practical calibration optimization can be obtained. Or example, in our case, in order to reach the same calibration noise as a common environmental clutter noise of 1-2 nT [1] only 60% coverage is required. Extrapolating the calibration time with the assumption of uniform calibration yields minimum of 4.2 minutes of rotating the sensor.

Of course, performing calibration in-situ cannot cover all orientations due to physical limitations of the human hand. In this case, data acquisition should be done systematically in order to approximate a uniform sampling.

## VII. Conclusion

A simple, low-cost, low-power method for performing a magnetic survey is presented. The method utilizes a low power consumption sensor and *ADC* with a calibration procedure that replaces the very expensive and cumbersome conventional scalar magnetometer. The setup of the survey system is quick and can be performed in-situ right before performing the magnetic survey, using the same hardware as a surveillance magnetometer. This approach allows an economical dual use of the same hardware as may be used in surveillance applications.

We present an *AI* based calibration algorithm and compare it to a conventional geometric calibration method. The *AI* scheme performs much better than the geometric one and allows us to reach 2 nT peak to peak value and less than $0.1 nT^2$ for the variance. In addition, we analyze the trade-off between calibration time and the output calibrated noise that can ultimately reach the theoretical sensor internal noise level or the practical environmental clutter.

The experimental results show the capability of using the calibration method for mapping purposes. The calibration procedure reduced the *RMS* noise from the order of $10^4$ nT to the order of 10 nT in complete 360º rotation in the natural earth magnetic field with variance of less than 0.1 nT. In many survey applications the sensor rotation is controlled manually to within ±10º resulting in a much lower PTP noise than in a complete 360º rotation. Hence, it is even more feasible to achieve magnetic mapping capabilities by using a low-cost, low-power fluxgate sensor.

## Acknowledgment

The authors would like to acknowledge Mr. Gal Zadok for contributing to the research and experiments in this work.